\documentclass{amsart}
\usepackage{latexsym,amsmath,amsfonts,amssymb,amsxtra}
\usepackage{color,latexsym,amsmath,amsthm,amsfonts,amssymb,xcolor,verbatim,lipsum,layout}

\usepackage{enumerate}

\usepackage{soul,cite,lineno,cancel}
\modulolinenumbers[2]

\scrollmode

\newtheorem{thm}{Theorem}[section]
\newtheorem{deff}[thm]{Definition}
\newtheorem{lem}[thm]{Lemma}

\newtheorem{rem}[thm]{Remark}
\newtheorem{prop}[thm]{Proposition}
\newtheorem{cor}[thm]{Corollary}

\newcommand{\ve}{\varepsilon}
\newcommand{\wt}{\widetilde}
\newcommand{\wh}{\widehat}
\newcommand{\vO}{\varOmega}

\newcommand{\ov}{\overline}

\def\N{{\mathbb N}}
\def\R{{\mathbb R}}

\def\mcL{{\mathcal L}}

\def\mcK{{\mathcal K}}

\def\mcN{\mathcal N}
\def\vL{{\varLambda}}
\def\mfA{{\mathfrak A}}
\def\mfB{{\mathfrak B}}
\def\mfC{{\mathfrak C}}
\def\mfD{{\mathfrak D}}
\def\mfF{{\mathfrak F}}

\def\mfG{{\mathfrak G}}
\def\mfM{\mathfrak M}
\def\mfN{\mathfrak N}
\def\mfZ{\mathfrak Z}

\def\mfW{\mathfrak W}
\def\mfN{\mathfrak N}

\def\mfP{\mathfrak P}

\def\emp{\emptyset}

\def\mcN{\mathcal N}
\def\smi{\setminus}
\def\bs{\boldsymbol}

\def\cd{\circledast}
\def\dd{\divideontimes}

\begin{document}

\title{Splitting of  Liftings in Product Spaces}
\author[K. Musia{\l}]{Kazimierz Musia\l}
\address{50-384 Wroc{\l}aw, Pl. Grunwaldzki 2/4, Poland}
\email{kazimierz.musial@math.uni.wroc.pl} 
\thanks{}

\date{\today}
\begin{abstract}
Let $(X, {\mathfrak A},P)$ and $(Y, {\mathfrak B},Q)$ be two probability spaces and $R$ be their skew product on the product $\sigma$-algebra ${\mathfrak A}\otimes\mfB$. Moreover, let
 $\{({\mathfrak A}_y,S_y)\colon y\in{Y}\}$ be a $Q$-disintegration of $R$ (if ${\mathfrak A}_y={\mathfrak A}$ for every $y\in{Y}$, then we have a regular conditional probability on ${\mathfrak A}$ with respect to $Q$) and let $\mfC$ be a sub-$\sigma$-algebra of ${\mathfrak A}\cap\bigcap_{y\in{Y}}{\mathfrak A}_y$. For  $f\in\mcL^{\infty}(R)$ I investigate the relationship between the $Y$-sections $[{\mathbb E}_{\mfC\otimes\mfB}(f)]^y$ of ${\mathbb E}_{\mfC\otimes\mfB}(f)$ (the conditional expectation of $f$ with respect to $\mfC\otimes\mfB$) and the conditional expectations of $f^y$ with respect $\mfC$ and $S_y$.

Moreover I prove the existence of a lifting $\pi$ on $\mcL^{\infty}(\wh{R})$ ($\wh{R}$ is the completion of $R$) and liftings $\sigma_y$ on $\mcL^{\infty}(\wh{S_y})$, $y\in Y$, such that
\begin{equation*}
[\pi(f)]^y= \sigma_y\Bigl([\pi(f)]^y\Bigr) \qquad\mbox{for all}
\quad y\in Y\quad\mbox{and}\quad f\in\mcL^{\infty}(\wh{R}).
\end{equation*}
As an application a characterization of stochastic processes possessing an equivalent measurable version is presented.
\end{abstract}


\date{\today}

\maketitle

\section{Preliminaries}
Throughout this paper $(X,\mfA,P)$ and $(Y,\mfB,Q)$ are fixed probability spaces and $\pi_X,\pi_Y$ denote the projections of $X\times{Y}$ onto $X$ and $Y$, respectively. $\mfA\times\mfB$ is the algebra on $X\times{Y}$ generated by $\pi_X^{-1}(\mfA)$ and $\pi_Y^{-1}(\mfB)$. $\mfA\otimes\mfB$ is the $\sigma$-algebra generated by $\mfA\times\mfB$. A measure $R$ on $(X\times{Y},\mfA\otimes\mfB)$ such that $\pi_X(R)=P$ and $\pi_Y(R)=Q$ is called a skew product of $P$ and $Q$.  We will use the notation $R=P\cd{Q}$.

If $(Z,\mfD,S)$ is a measure space, then  $\wh\mfD_S$ denote the completion od $\mfD$ with respect to $S$ and $\wh{S}$ denotes the completion of $S$. Moreover, $\mfD_0:=\{E\in\mfD:S(E)=0\}$ and $E^c:=Z\smi{E}$.  $S^*$ and $S_*$ are respectively the outer and the inner measures generated by $S$.
 $\mfA\wh\otimes_R\mfB$ is the completion of $\mfA\otimes\mfB$ with respect to $R$. If $\mfC$ is a sub-$\sigma$-algebra of $\mfD$ and $A\subset{Z}$, then $\mfC\ni{E}\supset{A}$ is called a $\mfC$-envelope of $A$ if $(S|_{\mfC})_*(E\smi{A})=0$.

$\vartheta(S)$ denotes the family of all lower densities on $(Z,\mfD,S)$ and $\vL(S)$ is the family of all liftings on the same space.

\begin{deff} \rm Let  $(X,\mfA)$ be a measurable space and $(Y,\mfB,Q)$ be a probability space. Let ${R}$ be a probability measure  on $\mfA\otimes\mfB$ such that $\pi_Y({R})=Q$. Assume that for every $y\in{Y}$ there exists a $\sigma$-algebra $\mfA_y\subset\mfA$ and a probability $Q_y$ on $\mfA_y$ with the following properties:
\begin{description}
\item [(Dis1)]
for each $A\in\mfA$ there exists  $N\in\mfB_0$ such that $A\in\mfA_y$ for all $y\notin{N}$ and the function $Y\smi{N}\ni y\longrightarrow S_y(A)$ is $\mfB$-measurable on $Y\smi{N}$;
\item[(Dis2)]
If $A\in\mfA$ and $B\in\mfB$, then
$$
\int_BS_y(A)\,dQ(y)={R}(A\times{B})\,.
$$
\end{description}
The family $\{(\mfA_y,S_y):y\in{Y}\}$ is called a  $Q$-disintegration of ${R}$.  If $\mfA_y=\mfA$ for all $y\in{Y}$, then $\{(\mfA,S_y):y\in{Y}\}$ is called a  regular conditional probability on $\mfA$ with respect to $Q$.

We say that $P$ is approximated by a family $\mfD\subset\mfA$ or $P$ is inner regular with respect to $\mfD$ if for each $A\in\mfA$ and $\ve>0$ there exists $D\in\mfD$ such that $D\subset{A}$ and $P(A\smi{D})<\ve$. \hfill$\Box$
\end{deff}
Disintegration not always exists. But if $P=\pi_X(R)$ is compact, then  it exists and possesses a nice approximation property.
\begin{prop}\label{p2} {\rm \cite[Theorem 3.5]{Pachl1}}
Let $(X,\mfA,P)$ and $Y,\mfB,Q)$ be probability spaces and $R=P\cd{Q}$. Assume that $Q$ is complete and $P$ is approximated by a countably compact lattice $\mcK\subset \mfA$ that is closed under countable intersections. Then there exists a $Q$-disintegration  $\{(\mfA_y,S_y):y\in{Y}\}$ of ${R}$ such that $\mfA_y\supset\mcK$ and $\mcK$ approximates $S_y$ for every $y\in{Y}$.
\end{prop}
Countably compact above means that if $\bigcap_{n=1}^{\infty}K_n=\emp$, then $\bigcap_{n=1}^mK_n=\emp$ for an $m\in\N$. Measures inner regular with respect to a countably compact collection of sets used to be called compact (Marczewski \cite{mar}).
\begin{rem} \rm Under the assumptions of Proposition \ref{p2} there exists a non-trivial $\sigma$-algebra $\sigma(\mcK)\subset \mfA\cap\bigcap_{y\in{Y}}\mfA_y$ approximating from inside all $\sigma$-algebras $\mfA_y$ and $\mfA$. In such a situation the $\sigma$-algebra $\mfC$ in Theorems \ref{T1}, \ref{T2} and \ref{T3} can be replaced by $\sigma(\mcK)$. In particular, \textbf{ Theorems \ref{T1}, \ref{T2} and \ref{T3} are valid in case of a compact $P$ and complete $Q$}.
\end{rem}
The following fact about skew products is known.
\begin{lem}\label{L2}
Let $(X,\mfA,P)$ and $Y,\mfB,Q)$ be probability spaces, $R=P\cd{Q}$ and  $\{(\mfA_y,S_y):y\in{Y}\}$ is  a  $Q$-disintegration of ${R}$.
If $f:X\times{Y}$ is a bounded  ${\mathfrak A}\otimes {\mathfrak
B}$-measurable function, then the function $y \in
{Y} \longmapsto \int_Xf^y(x)d\, S_y(x)$ is ${\mathfrak
B}$-measurable and the equality
$$
\int_{X \times Y}f(x,y)d\,R(x,y) = \int_Y\int_Xf^y(x)d\,S_y(x)dQ(y)
$$
holds true. In particular, if $E\in\mfA\otimes\mfB$, then $R(E)=\int_YS_y(E^y)\,dQ(y)$.

If $f$ is a bounded real valued ${\mathfrak A}\wh\otimes_R
{\mathfrak B}$-measurable function on $X \times Y$. Then
\begin{itemize}
\item[(a)]
the function $f^y$ is $\wh{\mathfrak A}_y$-measurable for $Q-a.e. \; y \in {Y}$;
 \item[(b)]
 if $f=0$ $\wh{R}$-a.e., then $f^y=0$ $\wh{S_y}$-a.e. for $Q-a.e. \; y \in {Y}$ ;
\item[(c)]
the function $y\longmapsto \int_Xf^y(x)d\wh S_y(x)$ is $\wh{\mathfrak B}$-measurable;
\item[(d)]
the equality
$$
\int_{X \times Y}f(x,y)d\wh R(x,y) =
 \int_Y\int_Xf^y(x)d\wh{S_y}(x)d\wh Q(y)
$$
holds true.
\end{itemize}
In particular, if $E\in{\mathfrak A}\wh\otimes_R{\mathfrak B}$,
then $\wh{R}(E)=\int_Y\wh{S_y}(E^y)\,d\wh{Q}(y)$.
\end{lem}
\section{ Conditional expectation in products }
We use the notation ${\mathbb E}_{\mfF}(f)$ for a version of the conditional expectation of $f$ with respect to $\mfF$ and $R$ and, ${\mathbb E}^y_{\mfG}(f^y)$ for a version of the conditional expectation of $f^y:=f(\cdot,y)$ with respect to $\mfG$ and $S_y$.

 Assume that $(X,\mfA,P)$  and $(Y,\mfB,Q)$ are probability spaces, $R=P\cd{Q}$ and $\{(\mfA,S_y):y\in{Y}\}$ is a regular conditional probability on $\mfA$ with respect to $\mfB$.  Let $\mfC$ be a sub-$\sigma$-algebra of $\mfA$ and $\mfM:=\mfC\otimes\mfB$. There is a natural question: what is the relation between ${\mathbb E}_{\mfM}(f)$ and the family $\{{\mathbb E}^y_{\mfC}(f^y)\colon y\in{Y}\}$? And more precisely: when $[{\mathbb E}^y_{\mfM}(f)]^y={\mathbb E}^y_{\mfC}(f^y)$ $S_y$-a.e. and for $Q$-almost all $y\in{Y}$ (the exceptional sets depend on $y$)? The first attempt may look the following: Let
\[
\mfW:=\{W\in\mfA\otimes\mfB \colon [{\mathbb E}_{\mfM}(\chi_W)]^y={\mathbb E}^y_{\mfC}([\chi_W]^y)\quad S_y-\mbox{a.e. and for }Q-\mbox{a.e. }y\in{Y}\}.
\]
$\mfW$ is a monotone class and remains only to prove that $\mfA\times\mfB\subset\mfW$. More precisely, it is enough to show that given $A\in\mfA$ one has $A\times{Y}\in\mfW$.

If $R$ is the direct product of $P$ and $Q$, then $S_y=P$ for $Q$-a.e. $y\in{Y}$ and the thesis is immediate (see \cite[Lemma 2.1]{mms2}). In the general case we come across the question of the $\mfM$-measurability of the function $(x,y)\longrightarrow {\mathbb E}^y_{\mfC}(\chi_A)(x)$.  If the function is $\mfM$-measurable, then the result follows from Lemma \ref{L2}(b). Unfortunately, it does not need to be so. Even more, we cannot a priori be sure that at least one such measurable family exists.

The subsequent theorem solves the problem in a more general situation.
\begin{thm}\label{T1} Assume that $(X,\mfA,P)$  and $(Y,\mfB,Q)$ are probability spaces, $R=P\cd{Q}$ and $\{(\mfA_y,S_y):y\in{Y}\}$ is a $Q$-disintegration of $R$. Let $\mfC$ be a sub-$\sigma$-algebra of $\mfA\cap\bigcap_{y\in{Y}}\mfA_y$  and $\mfM:=\mfC\otimes\mfB$.  If  $f \in {\mathcal L}^{\infty}(\wh{R})$ is arbitrary, then each version ${\mathbb E}_{\mfM}(f)$ of the conditional expectation of $f$ with respect to $\mfM$ is such that  $[{\mathbb E}^y_{\mfM}(f)]^y$ for  $Q$-almost all $y\in{Y}$  is a version of the conditional expectation ${\mathbb E}^y_{\mfC}(f^y)$ of $f^y$ with respect to $\mfC$ and $S_y$. Moreover, the function $(x,y)\longrightarrow {\mathbb E}^y_{\mfC}(f^y)(x)$ is then measurable with respect to $\mfM$ on a set $X\times(Y\smi{N_f})$, where $N_f\in\mfB_0$.
\end{thm}
\begin{proof}
Let us observe first that if $\mfN\subset \mfA\otimes\mfB$ is a $\sigma$-algebra and $h,g$ are two versions of ${\mathbb E}_{\mfN}(f)$, then $h=g \;\;R|_{\mfN}$-a.e. and then $h^y=g^y\;\;S_y$-a.e. for $Q$-almost all $y\in{Y}$ (Lemma \ref{L2}). Hence, during the inductive construction the subsequent conditional expectations ${\mathbb E}_{\mfM_{\alpha}}(f)$ will be chosen without any restrictions.

Without loss of generality, we may assume that  $\mfC$ is infinite (see \cite[Proposition 3.2]{smm}) and $f \in {\mathcal L}^{\infty}(R)$.  Let $\kappa$ be the smallest ordinal such that $\{K_{\alpha}: \alpha<\kappa\}$ is a (transfinite) sequence of elements of $\mfC$ such that $\mfC=\sigma\{K_{\alpha}: \alpha<\kappa\}$. Without loss of generality we may assume that $K_{\beta}\notin\mfC_{\beta}:=\sigma(\{K_{\alpha}\colon \alpha<\beta\})$, whenever $\beta<\kappa$ and  $\mfC_1:=\{\emp,X\}$. For each $\alpha<\kappa$ let
\[
\mfM_{\alpha}:=\mfC_{\alpha}\otimes\mfB\quad\mfC_{y0}:=\{C\in\mfC: S_y(C)=0\}\quad\mbox{and}\quad \mfC_{y\alpha}:=\sigma(\mfC_{\alpha}\cup\mfC_{y0})\quad\mbox{for }y\in{Y}.
\]
Then $\mfM_{\kappa}=\mfM$ and $\mfC_{\kappa}=\mfC$.

 $R_{\alpha}$ denotes the restriction of $R$ to $\mfM_{\alpha}$  and the restriction of $S_y$ to $\mfC_{y\alpha}$ by $S_{y\alpha}$.\\
 We will use the transfinite induction. Assume that for a $\gamma\leq\kappa$, for each $\alpha<\gamma$ and each $f\in\mcL^{\infty}(R)\subset\mcL^{\infty}(\wh{R})$ there exists a version ${\mathbb E}_{{\mfM}_{\alpha}}(f)$ of the conditional expectation of $f$ with respect to $\mfM_{\alpha}$ and $N_{f\alpha}\in\mfB_0$ such that  the equality ${\mathbb E}^y_{\mfC_{y\alpha}}(f^y):=[{\mathbb E}_{{\mfM}_{\alpha}}(f)]^y$ for  every $y\notin{N_{f\alpha}}$  defines a version of the conditional expectation of $f^y$ with respect to $\mfC_{y\alpha}$ and $S_y$.  \\
We begin with $\mfC_1=\{\emp,X\}$, $\mfC_{y1}=\sigma(\mfC_{y0})$ and $\mfM_1=\mfC_1\otimes\mfB$. If  $f \in {\mathcal L}^{\infty}(R)$  arbitrary, then - without loss of generality -  we may assume that the functions  $[{\mathbb E}_{\mfC_1\otimes\mfB}(f)]^y$  and ${\mathbb E}^y_{{\mathfrak C}_{y1}}(f^y)$ are constant for each $y\in{Y}$ separately.
If $B\in\mfB$, then (Lemma \ref{L2})
\begin{align*}
    \int_B \int_X \left[{\mathbb E}_{{\mathfrak M}_1}(f)\right]^y(x)\,d\wh{S_y}(x)\,d\wh{Q}(y)
        & =  \int_{X \times B} {\mathbb E}_{{\mathfrak M}_1}(f)(x,y)\,d\wh{R}(x,y)\\
        =\int_{X \times B} f(x,y)\, dR(x,y)
        & =  \int_B \int_{X} f^y(x)\, dS_y(x)\, dQ(y)\\
        = \int_B \int_{X} {\mathbb E}^y_{\mathfrak C_{y1}}(f^y)(x)\, dS_y(x)\, dQ(y)
&=\int_B \int_{X} {\mathbb E}^y_{\mathfrak C_{y1}}(f^y)(x)\, d\wh{S_y}(x)\, d\wh{Q}(y).
\end{align*}
This implies that  there exists $N\in {\mathfrak B}_0$ such that for any $y \in Y \smi{N}$ we have
\[
\int_X \left[{\mathbb E}_{{\mathfrak M}_1}(f)\right]^y(x)\, d\wh{S_y}(x) = \int_X {\mathbb E}^y_{\mathfrak C_{y1}}(f^y)(x)\,d\wh{S_y}(x)
\]
Since ${\mathbb E}^y_{\mathfrak C_{y1}}(f^y)$ and $\left[{\mathbb E}_{{\mathfrak M}_1}(f)\right]^y$ are constant, we have
 for all $y\in Y\setminus N$
$$
{\mathbb E}^y_{\mathfrak C_{y1}}(f^y)=\left[{\mathbb E}_{{\mathfrak M}_1}(f)\right]^y\quad\mbox{everywhere on } X\,.
$$
The further proof will be divided into three parts. \\
\textbf{(1)}$\bs{\gamma=\beta+1}$. \\
For simplicity, let $D:=K_{\beta}$. Then $\mfM_{\gamma}=\sigma(\mfM_{\beta}\cup\{D\times{Y}\})$ and $\mfC_{y\gamma}:=\sigma(\mfC_{y\beta}\cup\{D\})$. It is known that  $\mfC_{y\gamma}=\{(E\cap{D})\cup(F\cap{D^c})\colon E,F\in\mfC_{y\beta}\}$ and
$\mfM_{\gamma}=\{G\cap(D\times{Y})\cup{H}\cap(D^c\times{Y})\colon G,H\in\mfM_{\beta}\}$. Representations of single sets from $\mfM_{\gamma}$ and $\mfC_{y\gamma}$ are not unique. \\
If  $f$ is a fixed $\mfA\otimes\mfB$-measurable and bounded function, then there exist $\mfM_{\beta}$-measurable functions $f_1,f_2$ and $\mfC_{y\beta}$-measurable functions $f_{1y},f_{2y},\;y\in{Y}$, such that
\begin{align}
{\mathbb E}_{\mfM_{\gamma}}(f)&=f_1\chi_{D\times{Y}}+f_2\chi_{D^c\times{Y}}\quad R_{\gamma}-\mbox{a.e.}\label{e46}\\
{\mathbb E}^y_{\mfC_{y\gamma}}(f^y)&=f_{1y}\chi_D+f_{2y}\chi_{D^c}\quad S_{y\gamma}-\mbox{a.e.}\label{e47}
\end{align}
But
\begin{equation*}
{\mathbb E}_{\mfM_{\gamma}}(f)={\mathbb E}_{\mfM_{\gamma}}(f\chi_{D\times{Y}})+{\mathbb E}_{\mfM_{\gamma}}(f\chi_{D^c\times{Y}})\quad R_{\gamma}-\mbox{a.e.}
\end{equation*}
and
\begin{equation*}
{\mathbb E}^y_{\mfC_{y\gamma}}(f^y)={\mathbb E}^y_{\mfC_{y\gamma}}(f^y\chi_D)+{\mathbb E}^y_{\mfC_{y\gamma}}(f^y\chi_{D^c})\quad S_{y\gamma}-\mbox{a.e.}\,.
\end{equation*}
It follows from (\ref{e46}) and (\ref{e47}) that
\begin{align}
{\mathbb E}_{\mfM_{\gamma}}(f\chi_{D\times{Y}})&=\chi_{D\times{Y}}[f_1\chi_{D\times{Y}}+f_2\chi_{D^c\times{Y}}]= f_1\chi_{D\times{Y}}\quad R_{\gamma}-\mbox{a.e.}\label{e13}\\
{\mathbb E}^y_{\mfC_{y\gamma}}(f^y\chi_D)&=\chi_D[f_{1y}\chi_D+f_{2y}\chi_{D^c}]=f_{1y}\chi_D\quad S_{y\gamma}-\mbox{a.e.}\,.\label{e12}
\end{align}
Hence
\begin{align}
{\mathbb E}_{\mfM_{\beta}}(f\chi_{D\times{Y}})&=f_1{\mathbb E}_{\mfM_{\beta}}(\chi_{D\times{Y}})\quad R_{\beta}-\mbox{a.e.}\label{e9}\\
{\mathbb E}^y_{\mfC_{y\beta}}(f^y\chi_D)&=f_{1y}{\mathbb E}^y_{\mfC_{y\beta}}(\chi_D)\quad S_{y\beta}-\mbox{a.e.}\;y\in{Y}.\label{e10}
\end{align}
It follows from \eqref{e13} that if $B\in\mfB$, then
\begin{align*}
0&=\int_{X\times{B}}\left|[{\mathbb E}_{\mfM_{\gamma}}(f\chi_{D\times{Y}})]^y-[f_1\chi_{D\times{Y}}]^y\right|\,dR_{\gamma}\\
&=\int_B\int_X\left|[{\mathbb E}_{\mfM_{\gamma}}(f\chi_{D\times{Y}})]^y-f_1(\cdot,y)\chi_D\right|\,d S_{y\gamma}\,dQ(y)
\end{align*}
Hence there exists $N\in\mfB_0$ such that $\int_X|[{\mathbb E}_{\mfM_{\gamma}}(f\chi_{D\times{Y}})]^y-f_1(\cdot,y)\chi_D|\,d S_y=0$ if $y\notin{N}$. It follows that if $y\notin{N}$, then
\[
[{\mathbb E}_{\mfM_{\gamma}}(f\chi_{D\times{Y}})]^y=f_1(\cdot,y)\chi_D\quad S_{y\gamma}-a.e.
\]
Let $N_f^c:=\{y\in{Y}: [{\mathbb E}_{\mfM_{\gamma}}(f\chi_{D\times{Y}})]^y=f_1(\cdot,y)\chi_D\quad S_{y\gamma}-a.e.\}$. Then $N_f\in\mfB_0$.

By the assumption ${\mathbb E}^y_{\mfC_{y\beta}}(g^y)=[{\mathbb E}_{\mfM_{\beta}}(g)]^y$ everywhere for every  $g \in {\mathcal L}^{\infty}(R)$ for $Q$-almost all $y\in{Y}$ (the exceptional $Q$-null set depends on $g$) and so it follows from (\ref{e9}) and (\ref{e10}) that there exists $M_f\in\mfB_0$ such that  for each $y\notin{N_f}\cup{M_f}$ we have the equalities
\begin{align*}
f_{1y}{\mathbb E}^y_{\mfC_{y\beta}}(\chi_D)&\stackrel{\eqref{e10}}{=}{\mathbb E}^y_{\mfC_{y\beta}}(f^y\chi_D)=[{\mathbb E}_{\mfM_{\beta}}(f\chi_{D\times{Y}})]^y\\
&=f_1(\cdot,y)[{\mathbb E}_{\mfM_{\beta}}(\chi_{D\times{Y}})]^y=f_1(\cdot,y){\mathbb E}^y_{\mfC_{y\beta}}(f^y\chi_D)\quad S_{y\beta}-a.e.
\end{align*}
It follows that if $A_y:=\{x\in{X}: {\mathbb E}^y_{\mfC_{y\beta}}(\chi_D)(x)\neq0\}$, then $f_{1y}=f_1(\cdot,y)$ $S_{y\beta}$-a.e. on $A_y\in\mfC_{y\beta}$. That yields for each $y\notin{N_f}\cup{M_f}$ the equality
$$
[{\mathbb E}_{\mfM_{\gamma}}(f\chi_{D\times{Y}})]^y={\mathbb E}^y_{\mfC_{y\gamma}}(f^y\chi_D)\quad S_{y\gamma}-\mbox{a.e. on the set } A_y.
$$
Let us notice however that
$$
S_{y\gamma}(D\smi{A_y})=\int_{A_y^c}\chi_D\,dS_y=\int_{A_y^c}{\mathbb E}^y_{\mfC_{y\beta}}(\chi_D)\,dS_y=0
$$
and hence
$$
[{\mathbb E}_{\mfM_{\gamma}}(f\chi_{D\times{Y}})]^y={\mathbb E}^y_{\mfC_{y\gamma}}(f^y\chi_D)\quad S_{y\gamma}-\mbox{a.e. on }D.
$$
Taking into account the equations \eqref{e13} and \eqref{e12} for $y\notin{N_f\cup{M_f}}$, we obtain $S_{y\gamma}$-a.e. the required equality
$$
[{\mathbb E}_{\mfM_{\gamma}}(f\chi_{D\times{Y}})]^y= [f_1\chi_{D\times{Y}}]^y=f_1(\cdot,y)\chi_D=f_{1y}\chi_D={\mathbb E}^y_{\mfC_{y\gamma}}(f^y\chi_D)\,.
$$
In a similar way we obtain the equality
$$
[{\mathbb E}_{\mfM_{\gamma}}(f\chi_{D^c\times{Y}})]^y={\mathbb E}^y_{\mfC_{y\gamma}}(f^y\chi_{D^c})\quad S_{y\gamma}-\mbox{a.e.}
$$
Thus,
$[{\mathbb E}_{\mfM_{\gamma}}(f)]^y={\mathbb E}^y_{\mfC_{y\gamma}}(f^y)\quad S_{y\gamma}-\mbox{a.e.}$ for $Q$ almost every $y\in{Y}$ and so we may redefine ${\mathbb E}^y_{\mfC_{y\gamma}}(f^y)$ setting ${\mathbb E}^y_{\mfC_{y\gamma}}(f^y):=[{\mathbb E}_{\mfM_{\gamma}}(f)]^y$ everywhere on $X$ and for $Q$-almost every $y\in{Y}$.

 {\textbf{(2)} } $\bs{\gamma}$ \textbf{is of countable cofinality.}  Let a fixed sequence $(\gamma_n)_n$ be cofinal with $\gamma$. \\
  If $f\in{\mathcal L}^{\infty}(R)$, then the sequence $\{{\mathbb E}_{\mfM_{\gamma_n}}(f):n\in\N\}$ is a bounded martingale and so it is $R_{\gamma}$-a.e. and in $\mcL^1(R_{\gamma})$ convergent to ${\mathbb E}_{\mfM_{\gamma}}(f)$.
Similarly, for each $y\in{Y}$  the sequence $\{{\mathbb E}^y_{\mfC_{y\gamma_n}}(f^y):n\in\N\}$ is a bounded martingale. Hence, it is  $S_{y\gamma}$-a.e. and in $\mcL^1(S_{y\gamma})$ convergent to a version of ${\mathbb E}^y_{\mfC_{y\gamma}}(f^y)$. Moreover, by the construction, $[{\mathbb E}_{\mfM_{\gamma_n}}(f)]^y={\mathbb E}^y_{\mfC_{y\gamma_n}}(f^y)$ (everywhere on $X$ and) for a.e. $y\in{Y}$. Now we can verify the required equality for $\gamma$.
\begin{align*}
0&=\lim_n\int_{X\times{Y}}|{\mathbb E}_{\mfM_{\gamma_n}}(f)-{\mathbb E}_{\mfM_{\gamma}}(f)|\,dR_{\gamma}\\
&=\lim_n \int_Y\int_X\Bigl|[{\mathbb E}_{\mfM_{\gamma_n}}(f)]^y-[{\mathbb E}_{\mfM_{\gamma}}(f)]^y\Bigr|\,d{S_y}(x)\,d{Q}(y)\\
&= \lim_n \int_Y\int_X\Bigl|{\mathbb E}^y_{\mfC_{y\gamma_n}}(f^y)-[{\mathbb E}_{\mfM_{\gamma}}(f)]^y\Bigr|\,d{S_y}(x)\,d{Q}(y)\\
&\stackrel{Fatou}{\geq}\int_Y\lim_n\int_X\Bigl|{\mathbb E}^y_{\mfC_{y\gamma_n}}(f^y)-[{\mathbb E}_{\mfM_{\gamma}}(f)]^y\Bigr|\,d{S_y}(x)\,d{Q}(y)\\
&= \int_Y\int_X\Bigl|{\mathbb E}^y_{\mfC_{y\gamma}}(f^y)-[{\mathbb E}_{\mfM_{\gamma}}(f)]^y\Bigr|\,d{S_y}(x)\,d{Q}(y)\,.
\end{align*}
It follows that
\begin{equation*}
[{\mathbb E}_{\mfM_{\gamma}}(f)]^y=E^y_{\mfC_{y\gamma}}(f^y),\;\; S_{y\gamma}-a.e.,\;\mbox{and for }Q-\mbox{almost all }y\in{Y}\,.
\end{equation*}
Thus, one may set  ${\mathbb E}^y_{\mfC_{y\gamma}}(f^y):=[{\mathbb E}_{\mfM_{\gamma}}(f)]^y$  for $Q$-a.e. $y\in{Y}$.\\

 {\textbf{(3)} }  $\bs{\gamma}$ \textbf{is of uncountable cofinality.} Then $\mfC_{y\gamma}=\bigcup_{\alpha<\gamma}\mfC_{y\alpha}$. \\
If $f\in{\mathcal L}^{\infty}(R)$, then there exists a countably generated $\sigma$-algebra $\mfD\subset\mfC_{\gamma}$  such that $f$ and a version ${\mathbb E}_{\mfM_{\gamma}}(f)$ are measurable with respect to $\mfD\otimes\mfB$.  Then  for every $y\in{Y}$ the functions $f^y$ and $[{\mathbb E}_{\mfM_{\gamma}}(f)]^y$ are $\mfD$-measurable. Due to the uncountable cofinality of $\kappa$ we may assume that  $\mfD\subset\mfC_{\beta}$. In particular, $\mfD\subset\mfC_{y\beta}$, for every $y$. As a result, we may replace  all ${\mathbb E}_{\mfM_{\gamma}}(f)$ by ${\mathbb E}_{\mfM_{\beta}}(f)$. Similarly, for every $y\in{Y}$  the function ${\mathbb E}^y_{\mfC_{y\gamma}}(f^y)$  may be replaced by ${\mathbb E}^y_{\mfC_{y\beta}}(f^y)$. \\
In that way we have $[{\mathbb E}_{\mfM_{\gamma}}(f)]^y={\mathbb E}^y_{\mfC_{y\gamma}}(f^y)$ $S_{y\gamma}$-a.e. on $X$ and for $Q$-almost every $y\in{Y}$. We redefine ${\mathbb E}^y_{\mfC_{y\gamma}}(f^y)$ setting ${\mathbb E}^y_{\mfC_{y\gamma}}(f^y):=[{\mathbb E}_{\mfM_{\gamma}}(f)]^y$.

{\textbf{(4)}} As a result of the transfinite construction we have obtained the equality  $[{\mathbb E}_{\mfM_{\kappa}}(f)]^y={\mathbb E}^y_{\mfC_{\kappa}}(f^y)$ on $X$ for $Q$-a.e. $y\in{Y}$. If we have fixed at the very beginning a version ${\mathbb E}_{\mfM}(f)$, then we obtain the equality  ${\mathbb E}^y_{\mfC}(f^y):={\mathbb E}^y_{\mfC_{\kappa}}(f^y)=[{\mathbb E}_{\mfM_{\kappa}}(f)]^y= [{\mathbb E}_{\mfM}(f)]^y$ $S_y$-a.e. for $Q$-a.e. $y\in{Y}$.
The  measurability of the function $(x,y)\longrightarrow {\mathbb E}^y_{\mfC}(f^y)(x)$ is a direct consequence of the $R$-measurability of ${\mathbb E}_{\mfM}(f)$.
\end{proof}

As a corollary we obtain the following result:
\begin{thm}\label{T5}
Let $(X,\mfA,P)$ and $Y,\mfB,Q)$ be probability spaces and let $R=P\cd{Q}$. Assume that  $\{(\mfA,S_y):y\in{Y}\}$ is a regular conditional probability on $\mfA$ with respect to $Q$. Let $\mfC\subset\mfA$ be a $\sigma$-algebra and $\mfM=\mfC\otimes\mfB$. If  $f \in {\mathcal L}^{\infty}(\wh{R})$ is arbitrary, then each version ${\mathbb E}_{\mfM}(f)$ of the conditional expectation of $f$ with respect to $\mfM$ is such that  $[{\mathbb E}^y_{\mfM}(f)]^y$ for  $Q$-almost all $y\in{Y}$  defines a version of the conditional expectation ${\mathbb E}^y_{\mfC}(f^y)$ of $f^y$ with respect to $\mfC$ and $S_y$. Moreover, the function $(x,y)\longrightarrow {\mathbb E}^y_{\mfC}(f^y)(x)$ is then measurable with respect to $\mfM$ on a set $X\times(Y\smi{N_f})$, where $N_f\in\mfB_0$.
\end{thm}
\section{Liftings in product spaces}
We begin with an easy lemma from \cite[Lemma 2]{gvw}.
\begin{lem}\label{L3}
Let $(Z,\mfZ,T)$ be a probability space and let  $\mfC$ be a sub-$\sigma$-algebra of $\mfZ$ such that
$\mfZ_0\subset\mfC$. Moreover, let $\delta\in\vartheta(T|\mfC)$ be arbitrary
and let $M\in{\mfZ}\setminus\mfC$.
If $M_1\supset{M}$ and $M_2\supset{M^c}$ are $\mfC$-envelopes of $M$ and
$M^c$ respectively, then the formula
\begin{align*}
 &\widetilde\delta\Bigl[(G \cap M) \cup (H \cap M^c)\Bigr] :=  \\
   &  \Bigl[M \cap \delta\Bigl((G \cap M_1) \cup (H \cap
M^c_1)\Bigr)\Bigr] \cup \Bigl[M^c \cap\delta\Bigl((H \cap M_2) \cup (G
       \cap M^c_2)\Bigr)\Bigr]
\end{align*}
defines a density $\widetilde\delta\in\vartheta(T|\sigma(\mfC\cup\{M\}))$
that is an extension of $\delta$.
\end{lem}
It is worth to recall that the assumptions $\mfZ_0\subset\mfC$ (omitted in \cite{gvw}) is essential.
\begin{deff} \rm{ Let $(Z,\mfZ,T)$ be a probability space. A density
$\tau\in\vartheta(T)$ is called an {\em admissible density} (see \cite{mms2}) if it can be
constructed
with the help of the transfinite induction in the way described below.

\textbf{(A)}\quad
Let $\mfZ = \sigma\{M_{\alpha}:\alpha < \kappa\}$ be
numbered by the ordinals less than $\kappa$, where $\kappa$ is the first ordinal with this property and   for each $1 \leq \gamma < \kappa$ let  $\mfZ_{\gamma}=\sigma(\{M_{\alpha} :\alpha< \gamma \} \cup \mfZ_0)$.

\textbf{(B)}\quad
$\tau_0\in {\vartheta(T|\sigma(\mfZ_0))}$ is the only existing density on
$(Z,\sigma(\mfZ_0),T|\sigma(\mfZ_0))$, i.e.
\[
\tau_0(B)= \left\{ \begin{array}{ll}
\emp&if \ \  B \in {\mfZ_0} \\
&\\
Z &if \ \ B\notin {\mfZ_0}\ .
\end{array}\right.
\]
\textbf{(C)}\quad
If $\gamma<\kappa$ is a limit ordinal of uncountable cofinality,
then $\mfZ_{\gamma} = \bigcup_{\alpha < \gamma} \mfZ_{\alpha}$ and
we define $\tau_{\gamma}\in\vartheta(T|{\mfZ_{\gamma}})$ by setting
$$
\tau_{\gamma}(B) :\,=
\tau_{\alpha}(B)\quad\mbox{if}\quad B \in {\mfZ_{\alpha}}\quad\mbox{and}
\quad\alpha<\gamma\,.
$$
\textbf{(D)}\quad
Assume now that there exists an increasing sequence $(\gamma_n)$
of ordinals cofinal with $\gamma$.

For simplicity put $\tau_n :\,= \tau_{\gamma_n}$ and
$\mfZ_n :\,= \mfZ_{\gamma_n}$ for all $n \in {\bf N}$. Then $\mfZ_
{\gamma} = \sigma(\cup_{n \in {\bf N}}\mfZ_n)$ and we can define
$\tau_{\gamma}$ by setting
\begin{equation}\label{e20}
\tau_{\gamma}(B) := \bigcap_{k \in {\bf N}}\bigcup_{n \in {\bf N}}\bigcap_
{m \geq n}\tau_m(\{{\mathbb E}_{\mfZ_m}(\chi_B) > 1-1/k \})\quad\mbox{for}\quad
B \in {\mfZ_{\gamma}}\,.
\end{equation}
It follows from \cite[Lemma 1]{gvw} that $\tau_{\gamma} \in {\vartheta(T|\mfZ_{\gamma})}$ and $\tau_{\gamma} |\mfZ_n = \tau_n$
for each $n\in{\bf N}\,.$

\textbf{(E)}\quad
Let now $\gamma=\beta+1\,.$
To simplify the notations let $M : = M_{\beta}$.

It is well known, that
$$
\mfZ_{\gamma} = \{(G \cap M) \cup (H \cap M^c) : G, \, H \in \mfZ_{\beta}
\}\,.
$$
Let $M_1\supseteq M$ and $M_2\supseteq M^c$ be $\mfZ_{\beta}-$envelopes of
$M$ and $M^c$ respectively. We define a new density by
\begin{align*}
 &\tau_{\gamma}\Bigl((G \cap M) \cup (H \cap M^c)\Bigr) :=  \\
   &  \Bigl(M \cap \tau_{\beta}\Bigl((G \cap M_1) \cup (H \cap
M^c_1)\Bigr)\Bigr) \cup \Bigl(M^c \cap\tau_{\beta}\Bigl((H \cap M_2) \cup (G\cap M^c_2)\Bigr)\Bigr)
\end{align*}
for $G, H \in \mfZ_{\beta}\,.$
By Lemma \ref{L3}  $\tau_{\gamma} \in {\vartheta(T|\mfZ_{\gamma})}$ and
$\tau_{\gamma} |\mfZ_{\beta} = \tau_{\beta}\,.$

\textbf{(F)}\quad
We define $\tau\in\vartheta(T)$ by setting $\tau(E)=\tau_{\kappa}(E)$ for $E\in\mfZ$.} \hfill$\Box$
\end{deff}
\begin{rem}\rm
We might have assumed in the above definition that $M_{\gamma}\notin\mfZ_{\gamma}$ for all $\gamma$ but in Theorem \ref{T2} we need a definition free of such an assumption.
\end{rem}
\begin{deff} \rm Let $\tau\in\vartheta(\mu)$ be an admissible density. If $\rho\in\vL(\mu)$ is such that $\tau(A)\subset\rho(A)$ for every $A\in\mfZ$, then $\rho$ is called lifting  admissibly generated by $\tau$.
\end{deff}
Throughout, the collection of all admissible densities on
$(Z,\mfZ,T)$
will be denoted by $A\vartheta(T)$ and each $\tau\in A\vartheta(T)$
will be considered together with all elements involved into the above
construction without any additional remarks. In particular, for each limit ordinal $\gamma$ of countable cofinality, the cofinal sequence $\{\gamma_n:n\in\N\}$ will be fixed.
\hfill$\Box$

\begin{deff} \rm {Let $\{(X,\mfA_y,S_y)\colon y\in{Y}\}$ be a family of probability spaces and $(X,\mfA,P)$ be a probability space. Assume that \begin{align*}
\mfC:=&\sigma(\{M_{\alpha}\colon\alpha<\kappa\})\subset\mfA\cap\bigcap_{y\in{Y}}\mfA_y,\quad \mfC_1:=\{\emp,X\}\,,\mfC_{\gamma}:=\sigma(\{M_{\alpha}:\alpha<\gamma\})\\
\mfC_{y0}:=&\{C\in\mfC: S_y(C)=0\}\,,\quad\mbox{ and }
\quad \mfC_{y\gamma}:=\sigma(\mfC_{\gamma}\cup\mfC_{y0})\quad\mbox{for every }y\in{Y}.
\end{align*}
For each limit ordinal $\gamma$ of countable cofinality we fix one cofinal sequence $\{\gamma_n:n\in\N\}$. For each $y\in{Y}$ we are given an admissible density $\tau_y$ on $\mfC_y$, that is constructed along the transfinite sequence $\{M_{\alpha}\colon \alpha<\kappa\}$ with fixed  cofinal sequences $\{\gamma_n:n\in\N\}$ and for the transfinite sequence of $\sigma$-algebras $\mfC_{y\beta}$. It is important to notice that the construction of $\tau_y$ yields the inclusion $\tau_y(\mfC_{y\beta})\subset \mfC_{y\beta}$ for every $y\in{Y}$.  We call such a collection $\{\tau_y:y\in{Y}\}$ an equi-admissible family of densities determined by $\{M_{\alpha}\colon \alpha<\kappa\}$. }
\end{deff}
Moreover,  let
\[
\mfP_0:=\{W\in(\mfA\wh\otimes_R\mfB)_0\colon \forall \;y\in{Y}\;W^y\in\mfC\;\forall\;x\in{X}\;W_x\in\mfB\}\supset(\mfC\otimes\mfB)_0.
\]
\begin{thm}\label{T2}
Let $(X,\mfA,P)$  and $(Y,\mfB,Q)$ be probability spaces and $R$ be a skew product of $P$ and $Q$ with a disintegration $\{(\mfA_y,S_y):y\in{Y}\}$ with respect to $Q$. Assume that $\mfC\subset\mfA\cap\bigcap_{y\in{Y}}\mfA_y$ is such a $\sigma$-algebra that $P$ and each $S_y$ are inner regular with respect to $\mfC$. Assume moreover that $\{\tau_y:y\in{Y}\}$ is an equi-admissible family of densities determined by $\{M_{\alpha}\colon \alpha<\kappa\}$ generating $\mfC$.
Then, there
exists a density  $\varphi: \mfA\wh\otimes_R\mfB\to\sigma(\mfA\otimes\mfB\cup\mfP_0)$  such that for
each $F\in{\mfA}\wh\otimes_R{\mathfrak B}$
\begin{equation}\label{e18}
[\varphi(F)]^y = \tau_y([\varphi(F)]^y)\qquad\mbox{for  all }\;
y \in Y
\end{equation}
and
\begin{equation}\label{e19}
[\varphi(F)]_x\in\wh{\mathfrak B}\quad for\; all\; x\in{X}\,.
\end{equation}
\end{thm}
\begin{proof}  To prove the theorem, we will use the transfinite induction. \\
At the first step we  define $\varphi_1\in\vartheta\Bigl(\wh{R}|_{\sigma(\mfP_0)}\Bigr)$ by
$\varphi_1(E)=\emp$ if $\wh{R}(E)=0$ and $\varphi_1(E)=X\times{Y}$ otherwise.

To simplify further notation I denote each $\sigma$-algebra $\sigma({\mfC}_{\alpha}\otimes {\mathfrak B}\cup\mfP_0)$ by $\mfM_{\alpha}$.

Let us fix now $\gamma<\kappa$ and assume  that    for each $\alpha<\gamma$  there exists a density
\[
\varphi_{\alpha}: \sigma({\mfC}_{\alpha}\otimes\mfB\cup\mfP_0)\to \sigma({\mfC}_{\alpha}\otimes\mfB\cup\mfP_0)
\]
such that for each $F \in {\mfC}_{\alpha}\otimes {\mathfrak B}\cup\mfP_0$
\begin{equation}\label{e4}
[\varphi_{\alpha}(F)]^y\in\mfC_{y\alpha}\quad\mbox{and}\quad[\varphi_{\alpha}(F)]^y=\tau_y([\varphi_{\alpha}(F)]^y)\qquad\mbox{for every}\; y \in Y
\end{equation}
and
\begin{equation}\label{e3}
\varphi_{\beta}(F)=\varphi_{\alpha}(F)\quad \mbox{if }\alpha<\beta<\gamma\;\mbox{and\;}F\in\sigma({\mfC}_{\alpha}
\otimes {\mathfrak B}\cup\mfP_0).
\end{equation}
We will split now the proof into three parts:\\
\textbf{(A)} $\gamma=\beta+1$. \\
 Let $W_1\supset{M_{\beta}\times{Y}}$ and $W_2\supset{M_{\beta}^c\times{Y}}$ be respectively $\mfC_{\beta}\otimes\mfB$-envelopes of $M_{\beta}\times{Y}$ and $M_{\beta}^c\times{Y}$ with respect to $R$. Then, let $V_{1y}\supset{M_{\beta}}$ and $V_{2y}\supset{M_{\beta}^c}$ be respectively $\mfC_{y\beta}$-envelopes of $M_{\beta}$ and $M_{\beta}^c$ with respect to $S_y$. Without loss of generality, we may assume that
$$
V_{1y}\subset{W_1^y}\,,\quad V_{2y}\subset{W_2^y}\quad\mbox{for all }y\in{Y}\,.
$$
We define an extension $\ov\varphi_{\gamma}:\mfM_{\gamma}\to \mfM_{\gamma}$ of $\varphi_{\beta}$  to $\mfM_{\gamma}$ by the formula (see Lemma \ref{L3})
\begin{align*}
 &\ov\varphi_{\gamma}\Bigl[\bigl(G \cap (M_{\beta}\times{Y})\bigr) \cup \bigl(H \cap (M_{\beta}^c\times{Y})\bigr)\Bigr] :=  \\
   &  \Bigl[(M_{\beta}\times{Y}) \cap \varphi_{\beta}\Bigl((G \cap W_1) \cup (H \cap
W^c_1)\Bigr)\Bigr] \cup \Bigl[(M_{\beta}^c\times{Y}) \cap\varphi_{\beta}\Bigl((H \cap W_2) \cup (G
       \cap W^c_2)\Bigr)\Bigr]
\end{align*}
if $G,H\in\mfM_{\beta}$.

By the inductive assumption, if $y\in{Y}$, then
\begin{align*}
(\varphi_{\beta}[(G \cap W_1) \cup (H \cap W^c_1)])^y&=\tau_y[(\varphi_{\beta}[(G \cap W_1) \cup (H \cap W^c_1)])^y]\\
(\varphi_{\beta}[(H \cap W_2) \cup (G \cap W^c_2)])^y&=\tau_y[(\varphi_{\beta}[(H \cap W_2) \cup (G \cap W^c_2)])^y]
\end{align*}
The construction of the densities $\tau_y$ guarantees the equality
\begin{align}\label{e11}
 &\tau_y\Bigl[(A \cap M_{\beta}) \cup (B \cap M_{\beta}^c)\Bigr] :=  \\
   & \Bigl[M_{\beta} \cap \tau_y\Bigl((A \cap V_{1y}) \cup (B \cap
V^c_{1y})\Bigr)\Bigr] \cup \Bigl[M_{\beta}^c \cap\tau_y\Bigl((B \cap V_{2y}) \cup (A
       \cap V^c_{2y})\Bigr)\Bigr]\notag
\end{align}
if $A,B\in\mfC_{\beta}$.
It follows from  the basic properties of densities that for each $F\in\mfM_{\beta}$ there exists a set $N_1\in\mfB_0$ such that for every $y\notin{N_1}$ we have $\tau_y([\varphi_{\beta}(F)]^y)=\tau_y(F^y)$. Similarly, for each $F\in\mfM_{\gamma}$ there exists a set $N_2\in\mfB_0$ such that for every $y\notin{N_2}$ we have $\tau_y([\ov\varphi_{\gamma}(F)]^y)=\tau_y(F^y)$. Consequently, it follows from (\ref{e4}) and \eqref{e11} that  if $y\notin{N_1}\cup{N_2}$, then
\begin{align*}
 &\biggl[\ov\varphi_{\gamma}\Bigl[\bigl(G \cap (M_{\beta}\times{Y})\bigr) \cup \bigl(H \cap (M_{\beta}^c\times{Y})\bigr)\Bigr]\biggr]^y  \\
  &=\Bigl[M_{\beta} \cap \tau_y\Bigl((G^y \cap W_1^y) \cup (H^y \cap
(W_1^y)^c)\Bigr)\Bigr] \cup \Bigl[M_{\beta}^c \cap\tau_y\Bigl((H^y \cap W_2^y) \cup (G^y
       \cap (W_2^y)^c)\Bigr)\Bigr]\\
       & \supset \Bigl[M_{\beta} \cap \tau_y\Bigl((G^y \cap V_{1y}) \cup (H^y \cap
V^c_{1y})\Bigr)\Bigr] \cup \Bigl[M_{\beta}^c \cap\tau_y\Bigl((H^y \cap V_{2y}) \cup (G^y
       \cap V_{2y}^c)\Bigr)\Bigr]\\
       &=\tau_y\Bigl((G^y\cap{M_{\beta}})\cup (H^y\cap{M_{\beta}^c})\Bigr)= \tau_y\biggl[\Bigl[\bigl(G \cap (M_{\beta}\times{Y})\bigr) \cup \bigl(H \cap (M_{\beta}^c\times{Y})\bigr)\Bigr]^y\biggr]\\
         &=\tau_y\biggl(\biggl[\ov\varphi_{\gamma}\Bigl[\bigl(G \cap (M_{\beta}\times{Y})\bigr) \cup \bigl(H \cap (M_{\beta}^c\times{Y})\bigr)\Bigr]\biggr]^y\biggr).
\end{align*}
If $y\in{N_1}\cup{N_2}$, then $\biggl[\ov\varphi_{\gamma}\Bigl[\bigl(G \cap (M_{\beta}\times{Y})\bigr) \cup \bigl(H \cap (M_{\beta}^c\times{Y})\bigr)\Bigr]\biggr]^y\in\mfC_{y\gamma}$ and so we may  define
$\varphi_{\gamma}\in\vartheta(\wh{R}|_{{\mfM}_{\gamma}})$ setting for every $y\in{Y}$
 $[\varphi_{\gamma}(W)]^y=\tau_y([\ov\varphi_{\gamma}(W)]^y)$ for every $W\in{\mfM}_{\gamma}$. One can easily check that $\ov\varphi_{\gamma}|_{\mfM_{\beta}}=\ov\varphi_{\beta}$ and so $\varphi_{\gamma}|_{\mfM_{\beta}}=\varphi_{\beta}$.

\textbf{(B)} $\gamma$ is a limit ordinal of countable cofinality. \\
Assume that $(\gamma_m)_m$ is the increasing sequence of ordinals cofinal with $\gamma$.
Applying \eqref{e20} we define $\ov\varphi_{\gamma}\in\vartheta(\wh{R}|_{\mfM_{\gamma}})$ extending all $\varphi_{\gamma_m}$, setting for each $F\in\mfM_{\gamma}$
\begin{equation*}
\ov\varphi_{\gamma}(F):= \bigcap_{k \in {\bf N}}\bigcup_{n=1}^{\infty}
\bigcap_{m>n}\varphi_{\gamma_m}(\{{\mathbb E}_{{\mathfrak C}_{\gamma_m} \otimes
{\mathfrak B}}(\chi_F)
> 1 - 1/k\}).
\end{equation*}
According to \eqref{e20} we have also for each $A \in {\mfC_{\gamma}}$
\begin{equation*}
\tau_y(A):= \bigcap_{k \in {\bf N}}\bigcup_{n=1}^{\infty}
\bigcap_{m>n}\tau_y(\{{\mathbb E}^y_{{\mathfrak C}_{\gamma_m}}(\chi_A) > 1 -
1/k\})
\end{equation*}
 Then for $F \in \mfM_{\gamma}$ let $N_{m,F}\in{\mathfrak B}_0$ be such that (see
Theorem \ref{T1}) for each $y\notin{N_{m,F}}$ .
$$
[{\mathbb E}_{{\mathfrak C}_{\gamma_m} \otimes {\mathfrak
B}}(\chi_F)]^y={\mathbb E}^y_{{\mathfrak C}_{\gamma_m}}(\chi_{F^y})\;\;
S_y|{\mathfrak C}_{\gamma_m}-a.e.
$$
and (by the inductive assumption)
$$
[\varphi_{\gamma_m}(\{{\mathbb E}_{{\mathfrak C}_{\gamma_m} \otimes {\mathfrak B}}
       (\chi_F) >1 - 1/k\})]^y=\tau_y([\{{\mathbb E}_{{\mathfrak C}_{\gamma_m} \otimes {\mathfrak B}}
       (\chi_F) >1 - 1/k\}]^y.
$$
If $ N_F:= \bigcup_{m \in {\bf N}} N_{m,F}$, then $N_F \in
{\mathfrak B}_0$ and it follows from Theorem \ref{T1}, that
for all $y \in Y \setminus N_F$
\begin{align*}
    [\ov\varphi_{\gamma}(F)]^y & =  \bigcap_{k \in {\bf N}}\bigcup_{n \in {\bf N}}
      \bigcap_{m \geq n}[\varphi_{\gamma_m}(\{{\mathbb E}_{{\mathfrak C}_{\gamma_m} \otimes {\mathfrak B}}
       (\chi_F) >1 - 1/k\})]^y \\
        & =  \bigcap_{k \in {\bf N}}\bigcup_{n \in {\bf N}}\bigcap_{m
              \geq n}\tau_y([\{{\mathbb E}_{{\mathfrak C}_{\gamma_m} \otimes {\mathfrak B}}
              (\chi_F) > 1 - 1/k\}]^y) \\
        & =  \bigcap_{k \in {\bf N}}\bigcup_{n \in {\bf N}}\bigcap_{m
              \geq n}\tau_y(\{[{\mathbb E}_{{\mathfrak C}_{\gamma_m} \otimes {\mathfrak B}}
              (\chi_F)]^y > 1 - 1/k\}) \\
        & =  \bigcap_{k \in {\bf N}}\bigcup_{n \in{\bf N}}\bigcap_{m \geq n}
              \tau_y(\{{\mathbb E}^y_{{\mathfrak C}_{\gamma_m}}(\chi_{F^y}) > 1 - 1/k\}) =
              \tau_y(F^y).
\end{align*}
It follows that $\tau_y([\ov\varphi_{\gamma}(F)]^y)=[\ov\varphi_{\gamma}(F)]^y$,
for all $y \notin M_F\in{\mathfrak B}_0$, where $N_F\subset{M_F}$.  We define $\varphi_{\gamma}$ by setting $[\varphi_{\gamma}(F)]^y :=\tau_y([\ov\varphi_{\gamma}(F)]^y)$ for every $y\in{Y}$, because if $y\in{M_F}$, then $[\ov\varphi_{\gamma}(F)]^y\in\mfC_{y\gamma}$.

\textbf{C)} Assume now that $\gamma$ is of uncountable cofinality. \\
In such a case $\mfC_{\gamma}=\bigcup_{\alpha<\gamma}\mfC_{\alpha}$ and for each $F\in\mfC_{\gamma}$ there exists $\alpha<\gamma$ with $F\in\mfC_{\alpha}$. We simply set $\varphi(F)=\varphi_{\alpha}(F)$.\\
\textbf{D)} When $\gamma=\kappa$, then we get a density $\varphi\in\vartheta(\wh{R}|_{\mfM_{\kappa}})=\vartheta(\wh{R}|_{\sigma(\mfC\otimes\mfB\cup\mfP_0)}$ and the inductive part is finished.

\textbf{D)} Now we are going to extend it to a density $\varphi:\mfA\wh\otimes_R\mfB\to \mfA\otimes\mfB\cup\mfP_0$.
Let
\[
\mfW:=\{W\in\mfA\wh\otimes_R\mfB\colon \exists\;V\in\mfC\otimes\mfB\;\&\;\wh{R}(W\triangle{V})=0\}
\]
We define a density  $\varphi:\mfA\wh\otimes_R\mfB\to \mfC\otimes\mfB\cup\mfP_0$ setting $\varphi(W):=\varphi(V)$.
In that way the equality (\ref{e18}) holds true.

To see that (\ref{e19}) is also valid it is enough to observe that if $W$ and $V$ are as above, then $[\varphi(V)]_x\in\wh{\mfB}$ for every $x\in{X}$, due to the corresponding property of $\mfP_0$.
\end{proof}
\begin{prop}\label{p3}
Under the assumptions of Theorem \ref{T2}
there exists $\psi\in\vartheta( \wh{R})$ satisfying for all $E\in{\mathfrak A}\wh\otimes_R {\mathfrak B}$ the following conditions:
\begin {itemize}
\item[(1)]
$ \wh{S_y}([\psi(E)]^y\cup[\psi(E^c)]^y)=1 \qquad\mbox{for
all}\quad y\in Y\,;$
\item[(2)]
$[\psi(E)]^y= \tau_y\Bigl([\psi(E)]^y\Bigr) \qquad\mbox{for all}
\quad  y\in Y\,;$
\item[(3)]
$[\psi(E)]_x \quad\mbox{is}\quad \wh{Q}-\mbox{measurable for
all}\quad x\in X\,.$
\end{itemize}
\end{prop}
\begin{proof}  For each $y\in{Y}$ the measure $S_y$ is inner regular with respect to $\mfC$ and consequently the density $\tau_y\in\vartheta(S_y|_{\mfC})$ can be uniquely extended to $\wh{S_y}$. We denote the extension also by $\tau_y$. Take $\varphi\in\vartheta( \wh{R})$ satisfying the thesis
of Theorem \ref{T2} and let
\begin{equation*}
\Phi:\,=\biggl\{\overline\varphi\in\vartheta( \wh{R}):\forall y\in Y
\,\forall E\in{\mathfrak A}\wh\otimes_R{\mathfrak
B}\;\varphi(E)\subseteq\overline\varphi(E)\;\&\;[\overline\varphi(E)]^y\subseteq
\tau_y([\overline\varphi(E)]^y)\biggr\}.
\end{equation*}
We consider $\Phi$ with inclusion as the partial order:
$\overline\varphi_1\leq\overline\varphi_2$ if
$\overline\varphi_1(E)\subseteq\overline\varphi_2(E)$ for each
$E\in{\mathfrak A}\wh\otimes_R{\mathfrak B}$. As in the proof of  \cite[Theorem 3.5]{smm} one can show that the maximal element of $\Phi$ satisfies the thesis of Proposition \ref{p3}.
\end{proof}

\begin{thm}\label{T3}
Let $(X,\mfA,P)$  and $(Y,\mfB,Q)$ be probability spaces. Assume that $R=P\cd{Q}$ has a disintegration $\{(\mfA_y,S_y):y\in{Y}\}$ with respect to $Q$ and $\mfC\subset\mfA\cap\bigcap_{y\in{Y}}\mfA_y$ is such a $\sigma$-algebra that $P$ and each $S_y$ are inner regular with respect to $\mfC$. Let $\{\tau_y\in\vartheta(\wh{S_y}):y\in{Y}\}$ be the family from Proposition \ref{p3}.  Then there exist liftings $\sigma_y\in\vL(\wh{S_y})$ that are admissibly generated by $\tau_y$ for all $y\in Y$ and there exists $\pi\in {\vL(\wh R)}$ such that the following condition is satisfied:
\begin{equation*}
[\pi(E)]^y= \sigma_y\Bigl([\pi(E)]^y\Bigr) \qquad\mbox{for all}
\quad y\in Y\quad\mbox{and}\quad E\in {\mathfrak
A}\wh\otimes_R{\mathfrak B}.
\end{equation*}
Equivalently,  for each $f \in {\mathcal L}^{\infty}(\wh R)$ and each
$y\in Y$
$$
[\pi(f)]^y = \sigma_y\Bigl([\pi(f)]^y\Bigr).
$$
\end{thm}
\begin{proof}   According to Proposition \ref{p3} there exists
 $\psi \in{\vartheta(\wh R)}$ such that for all $E \in {\mathfrak A}\wh\otimes_R{\mathfrak B}$
\begin{equation*}
[\psi(E)]^y\subset\tau_y\Bigl([\psi(E) ]^y\Bigr)\qquad\hbox{ for all
$y\in Y\,$}\,,
\end{equation*}
and
\begin{equation}\label{7b}
\wh{S_y}\Bigl([\psi(E)]^y\cup[\psi(E^c)]^y \Bigr)=1
\qquad\hbox{for all}\quad y\in Y\;.
\end{equation}

We take now for each $y\in Y$ a lifting $\sigma_y\in\vL(\wh{S_y})$
such that $\tau_y\subseteq\sigma_y$ and define
$\pi\in\vartheta(\wh R)$ by setting for each $E\in {\mathfrak
A}\wh\otimes_R{\mathfrak B}$ and each $y\in Y$
\begin{equation}\label{19}
[\pi(E)]^y= \sigma_y\Bigl([\psi(E)]^y\Bigr)\,.
\end{equation}
Since $\psi(E)\subseteq\pi(E)$ for all $E\in{\mathfrak
A}\wh\otimes_R{\mathfrak B}$ we get $\wh{R}$-measurability of
$\pi(E)$ and $\pi(E)\stackrel{\wh{R}}{=}E$. In order to prove that
$\pi$ is a lifting it suffices to show that we have always  $
\pi(E^c)=[\pi(E)]^c$. But this is a consequence of (\ref{7b}) and
(\ref{19}) as we get for each $y$ the equality
$[\pi(E^c)]^y=\Bigl([\pi(E)]^y\Bigr)^c\,$. This proves that
$\pi\in\vL(\wh R)\,.$
\end{proof}
\begin{cor}\label{c1}
Under the assumptions of Theorem \ref{T3} if $E\in\mfA\wh\otimes_R\mfB$ and $\wt{E}$ is defined by $\wt{E}^y:=\sigma_y(E^y)\,, y\in{Y}$, then $\wt{E}\in\mfA\wh\otimes_R\mfB$.
\end{cor}
\begin{proof}
According to Theorem \ref{T3} $[\pi(E)]^y= \sigma_y\Bigl([\pi(E)]^y\Bigr)$, for every $y\in{Y}$. In virtue of Lemma \ref{L2} there exists $N\in\mfB_0$ such that $[\pi(E)]^y= \sigma_y(E^y)=\wt{E}^y$ if $y\notin{N}$. As a result $\wt{E}\cap(X\times{Y\smi{N}})=\pi(E)\cap(X\times{Y\smi{N}})\in\mfA\wh\otimes_R\mfB$.
\end{proof}
\begin{deff} \rm
Let $\mcN:=\{E\subset{X\times{Y}}\colon \exists\;N\in\mfB_0\;\forall\;y\notin{N}\;\wh{S_y}(E^y)=0\}$. Elements of $\mcN$ are called $R$-nil sets and form a $\sigma$-ideal (compare with \cite[\S1,II]{Haupt}). Then, let $\mfA\dd_R\mfB:=\{W\triangle{N}\colon W\in\mfA\otimes\mfB\;\;\&\; N\in\mcN\}$ and $R_{\dd}(W\triangle{N}):=R(W)$ be the nil extensions of $\mfA\otimes\mfB$ and $R$.
$R_{\dd}$ is a complete measure extending $\wh{R}$. If $\pi_1\in\vL(\wh{R})$, then the formula $\pi_2(W\triangle{N}):=\pi_1(W)$ defines a lifting for $R_{\dd}$.
\end{deff}

The following theorem is a direct consequence of Theorem \ref{T3} and it is a strong generalization of \cite[Theorem 3.6]{smm}, where the $\sigma$-algebra $\mfA$ was assumed to be separable in the  Frechet-Nikod\'{y}m pseudometric.
\begin{thm}\label{T4}
Let $(X,\mfA,P)$  and $(Y,\mfB,Q)$ be probability spaces and  $R$ be a skew product of $P$ and $Q$.  Assume that   $\{(\mfA,S_y):y\in{Y}\}$ is a regular conditional probability on $\mfA$ with respect to $Q$  and  $\{\tau_y\in\vartheta(S_y):y\in{Y}\}$ is an equi-admissible family of densities. Then for each $y\in{Y}$ there exists $\sigma_y\in\vL(\wh{S_y})$ admissibly generated by $\tau_y$ and there exists a lifting $\pi:\mfA\dd_R\mfB\to \mfA\wh\otimes_R\mfB$ such that the following condition is satisfied:
\begin{equation*}
[\pi(E)]^y= \sigma_y\Bigl([\pi(E)]^y\Bigr)\;\;\mbox{and }E\triangle\pi(E)\in\mcN  \quad\mbox{for all}
\;\; y\in Y\;\;\mbox{and}\;\; E\in {\mathfrak
A}\dd_R{\mathfrak B}.
\end{equation*}
Equivalently,  for each $f \in {\mathcal L}^{\infty}(R_{\dd})$ and each
$y\in Y$
$$
[\pi(f)]^y = \sigma_y\Bigl([\pi(f)]^y\Bigr)\quad\mbox{and}\;\{(x,y):\pi(f)(x,y)\neq{f(x,y)}\}\in\mcN.
$$
\end{thm}
Theorem \ref{T4} is also a partial generalization of the main results from \cite{mms1} and \cite{mms2}. In \cite{mms2} the measure $R$ was a direct product of $P$ and $Q$ whereas in \cite{mms1} the measure $R$ was assumed to be absolutely continuous with respect to the product $P\times{Q}$. With stronger assumptions one could achieve stronger results. The main one was - so called - rectangle formula. In case of the direct product it looks as $\pi(A\times{B})=\rho(A)\times\sigma(B)$, where liftings $\rho$ and $\sigma$ are properly chosen. If $R\ll{P\times{Q}}$, then $\pi(A\times{B})=\bigcup_{y\in\rho(B)}(\sigma_y(A)\times\{y\})$, where again the liftings are properly chosen.
\section{Measurability of stochastic processes}
Let $(X,\mfA,P),\,(Y,\mfB,Q),\;R,\;R_{\dd}$ and the regular conditional probability $\{(S_y,\mfA),\;y\in{Y}\} $ be as in Theorem \ref{T4}. Moreover, let  $\{\xi_y:y\in{Y}\}$ and  $\{\theta_y:y\in{Y}\}$ be stochastic processes defined on $(X,\mfA,P)$. We say that the processes are equivalent if $S_y\{x\in{X}:\xi_y(x)\neq\theta_y(x)\}=0$ for every $y\in{Y}$. The process $\{\theta_y:y\in{Y}\}$ is called measurable if the function $(x,y)\longrightarrow\theta_y(x)$ is $\wh{R}$-measurable.

There are known several examples of stochastic processes without equivalent measurable version (cf. \cite{Cohn, Dudley, HJ}).
The following theorem describes the stochastic processes possessing such a measurable modification. I consider only real valued processes.
\begin{prop}
A  stochastic process $\{\xi_y:y\in{Y}\}$  possesses an equivalent measurable version if and only if  the process is $R_{\dd}$-measurable. Every modification of a measurable process with the help of  admissibly generated liftings is a measurable process.
\end{prop}
\begin{proof}
Assume at the beginning that the process defined by $\varXi(x,y):=\xi_y(x)$ is uniformly bounded.

Let $\varTheta(x,y):=\theta_y(x)$ be measurable and equivalent to $\{\xi_y:y\in{Y}\}$. Define the set $N$ by the equality $N:=\{(x,y)\in{X\times{Y}}\colon \theta_y(x)\neq\xi_y(x)\}$.  By the assumption $P(N^y)=0$ for every $y\in{Y}$. Thus $N\in\mcN$ what means that $\{\xi_y:y\in{Y}\}$ is $R_{\dd}$-measurable. If $B\subset\R$ is a Borel set, then $\varXi^{-1}(B)=\varXi^{-1}(B)\cap{N^c}\cup\varXi^{-1}(B)\cap{N}=\varTheta^{-1}(B)\cap{N^c}\cup\varXi^{-1}(B)\cap{N}$ and consequently $\varXi^{-1}(B)\triangle\varTheta^{-1}(B)\subset{N}\in\mcN$.

And the reverse implication. If  the function $\varXi$ is $R_{\dd}$-measurable, then we apply Theorem \ref{T4}. It follows from the elementary properties of lifting that there exists a set $N_1\in\mfB_0$ such that $[\pi(\varXi)]^y=\sigma_y(\varXi^y)$ $\wh{S_y}$-a.e. and for every $y\notin{N_1}$.  We define the process $\{\theta_y:y\in{Y}\}$ setting $\theta_y:=\sigma_y(\Bigl([\pi(\varXi)]^y\Bigr)$ if $y\notin{N_1}$ and $\theta_y:=\xi_y$ if $y\in{N_1}$. In virtue of Theorem \ref{T4} the function $\pi(\varXi)$ is $\wh{R}$-measurable.

If the process $\varXi$ is arbitrary, then one can decompose the space $X\times{Y}$ into at most countably many pairwise disjoint sets $\vO_n\in\mfA\dd_R\mfB$ such that $\varXi$ is bounded on each $\vO_n,\,n\in\N$. Let for each $n\in\N$ $\varXi_n(x,y):=\varXi(x,y)\,\chi_{\vO_n}(x,y)$. Then, we  define a measurable version $\varTheta_n$ of $\varXi_n$ setting  $\varTheta_n(x,y):=\sigma_y(\Bigl([\pi(\varXi_n)]^y\Bigr)$ if $y\notin{N_n}$ and $\varTheta_n(\cdot,y):=\xi_y(\cdot)\,\chi_{[\pi(\vO_n)]^y}(\cdot)$ if $y\in{N_n}$, where $N_n\in\mfB_0$ and setting $[\pi(\varTheta_n)]^y=\sigma_y(\varTheta_n^y),\;\wh{S_y}$-a.e.  for every $y\notin{N_n}$. Since the supports of the functions $\varTheta_n$ are pairwise disjoint, we may define the final measurable process $\varTheta$ setting
\[
\varTheta(x,y) = \left\{ \begin{array}{ll}
\varTheta_n(x,y)&{\rm\; if} \ \  (x,y)\in\pi(\vO_n)\cap(X\times{N^c_n})\\
&\\
\xi_y(x)&{\rm\; if} \ \  (x,y)\in\bigcup_n(X\times{N_n})\cup\bigcap_n\pi(\vO_n^c)\ .   \\
\end{array}\right.
\]
\end{proof}

\end{document}